\documentclass{article}
\usepackage{latexsym}
\usepackage{amssymb}
\usepackage{amsmath}
\usepackage{color}
\usepackage{graphicx}
\allowdisplaybreaks

\renewcommand\thesection{\arabic{section}}
\renewcommand\theequation{\thesection.\arabic{equation}}

\newtheorem{theorem}{\color{black} Theorem}[section]
\newtheorem{lemma}{\color{black} Lemma}[section]

\newtheorem{corollary}{\color{black} Corollary}[section]

\begin{document}
\large
\title{\bf Persistence of lower dimensional invariant tori on sub-manifolds
in Hamiltonian systems}

\author{
Zhenxin Liu\thanks{{\it E-mail address:} zxliu@email.jlu.edu.cn
(Zhenxin Liu).}
\\{\it School of Mathematics, Jilin University,
Changchun 130012, P. R. China} }
\date{}
\maketitle {\noindent{\bf Abstract}}

Chow, Li and Yi  in \cite{cl} proved that the majority of the
unperturbed tori {\it on sub-manifolds} will persist for standard
Hamiltonian systems. Motivated by their work, in this paper, we
study the persistence and tangent frequencies preservation of
lower dimensional invariant tori on smooth sub-manifolds for real
analytic, nearly integrable Hamiltonian systems. The surviving
tori might be elliptic, hyperbolic, or of mixed type.\\
 {\it keywords}: Hamiltonian system; lower dimensional invariant
 tori;
persistence on sub-manifolds; KAM theorem

\indent

\vskip 5mm
\section{Introduction}

 The persistence of quasi-periodic solutions or invariant tori to
integrable Hamiltonian systems had puzzled scientists for long up
to the appearance of the celebrated KAM \cite{kol,arn,mo} theory,
which affirmed that the majority of invariant tori persist under
small perturbations.

Later Melnikov \cite{mel} formulated a  KAM type persistence
result for elliptic lower dimensional tori of integrable
Hamiltonian systems under  so-called Melnikov's non-resonance
condition. More precisely, for a system with the following
Hamiltonian
$$H=N+P=\sum_{j=1}^{n}\omega_{j}y_{j}+\frac{1}{2}\sum_{j=1}^{m}\Omega_{j}
(u_{j}^{2}+v_{j}^{2})+P,$$ Melnikov announced that the majority of
invariant tori survive the small perturbations  under the
following conditions
$$|\langle k,\omega\rangle+\langle l,\Omega\rangle|>\frac{\gamma}
{|k|^{\tau}},\ |l|\leq 2$$
for $k\in Z^{n},\ l\in Z^{m},|k|+|l|\neq 0.$ The complete proof of
his result was later carried out by Eliasson, Kuksin, and
P\"{o}schel \cite{eli, kuk, p}. In fact, Moser \cite{mos} had
already noted the persistence of elliptic lower dimensional tori.
He proved that the existence of the quasi-periodic solutions when
the tori admit 2-dimensional elliptic equilibrium point. However
the way he used can not be applied to higher dimension, because he
requested the tangent frequency be fixed. Later, Eliasson
\cite{eli} removed the restriction successfully by letting the
frequency suffer small perturbations and later P\"{o}schel
\cite{p} simplified the proof of Eliasson \cite{eli}.

For Hamiltonian
$$H=N+P=\langle\omega_{0},y\rangle+\frac{1}{2}\langle u,Mu\rangle+P(x,y,u),$$
where $(x,y,u)\in T^{n}\times R^{n}\times R^{2m}$, Moser
\cite{mos} obtained the persistence of hyperbolic invariant tori
when $\omega_0\in R^n$ is a fixed Diophantine toral frequency and
the eigenvalues of $JM$ ($J$ being the standard symplectic matrix
in $R^{2m}$) are real and distinct. Graff \cite{gra} generalized
Moser's result by allowing multiple eigenvalues of $JM$. And the
proof of Graff's result was later given by Zehnder \cite{zeh}, who
used implicit function techniques. More recently, Li and Yi
\cite{ly2} generalized the results of Graff and Zehnder on the
persistence of hyperbolic invariant tori in  Hamiltonian systems
by allowing the degeneracy of the unperturbed Hamiltonians and
they obtain the preservation of part or full components of tangent
frequencies. They adopted the Fourier series expansion for normal
form $N$, which is a new technique.

Recently, Chow, Li and Yi \cite{cl} proved that the majority of
the unperturbed tori {\it on sub-manifolds} will persist under a
non-degenerate condition of R\"{u}ssmann type for standard
Hamiltonian systems. Motivated by their work, in this paper, we
shall show that lower dimensional tori also survive small
perturbations on sub-manifolds under some assumptions. The
surviving tori might be elliptic, hyperbolic, or of mixed type.

We consider a real analytic family of Hamiltonian systems of the
following form
\begin{equation}
H(x,y,u)=N(y,u)+P(x,y,u),
\end{equation}
where $(x,y,u)$ lies in a complex neighborhood $\{(x,y,u):|{\rm
Im}x|\leq r,{\rm dist}(y,G)\leq s,|u|\leq s\}$ of $T^{n}\times
G\times \{0\}\subset T^{n}\times R^{n}\times R^{2m},\ G\subset
R^{n}(n\geq 2)$ is a bounded closed region and $P$ is small.
Besides these, we also assume that

\begin{itemize}
\item[{\bf A0)}] $N_{u}(y,0)=0,\ {\rm det}N_{uu}(y,0)\neq 0.$
\end{itemize}

To prove the persistence of lower dimensional invariant tori of
system (1.1), we first consider the following parameter-dependent,
real analytic Hamiltonian system
\begin{equation}
H=e(\lambda)+\langle\omega (\lambda),y\rangle+\frac{1}{2}\langle
A(\lambda)y,y\rangle+\frac{1}{2} \langle
M(\lambda)u,u\rangle+P(x,y,u,\lambda),
\end{equation}
where $(x,y,u)$ lies in a complex neighborhood $D(r,s)=\{
(x,y,u):|{\rm Im}x|\leq r,|y|\leq s, |u|\leq s\}$ of $T^{n}\times
\{0\}\times \{0\}\subset T^{n}\times R^{n}\times R^{2m}$,
$\lambda$ is a parameter lying in a bounded closed region $\Lambda
\subset R^{n_{0}}$ and $M(\lambda)$ is nonsingular on $\Lambda$.
In the above, all $\lambda$ dependency are of class $C^{l_{0}}$for
some $l_{0}\geq n$.

We assume the following conditions:

\begin{itemize}
\item[{\bf A1)}]  rank$\{{{\partial^{\alpha}\omega}\over {\partial
\lambda ^{\alpha}}}:|\alpha |\leq n-1\}=n$  for all $\lambda \in
\Lambda.$

\item[{\bf A2)}]  For $(k,l)$ satisfying $0<|k|\leq
K=\frac{4n}{\sigma}{\rm max}_{1\leq i\leq 2m,0\leq r\leq n-1}
|\partial_{\lambda}^{r}\Omega_{i}|$ and $|l|\leq 2,$ the following
holds:
$${\rm meas}\{\lambda:|i\langle k,\omega(\lambda)\rangle+\langle l,
\Omega(\lambda)\rangle|=0\}=0,$$ where
$\Omega=(\Omega_{1},\cdots,\Omega_{2m})^{\top}$, $\Omega_{j},\
j=1,\cdots,2m$ are eigenvalues of $JM$, $\sigma$ is a constant
which will be determined later and ``meas" denotes the Lebesgue
measure in $R^{n_{0}}.$

\item[{\bf A3)}]  rank$A(\lambda)=d$  on $\Lambda$, and, there is
a smoothly varying, nonsingular, $d\times d$ principal minor
$\tilde{A} (\lambda)$ of $A(\lambda)$.
\end{itemize}

Denote $i_{1},i_{2},\cdots,i_{d}$ as the row indices (in the
natural order) of $\tilde{A} (\lambda)$ in $A(\lambda)$.

Our main result states as follows.

\begin{theorem}Consider (1.2).

1) Assume A1), A2), A3) hold and let $\tau>n(n-1)-1$ be fixed.
Then for a given $\gamma$ there exists an
  $\epsilon=\epsilon(r,s,l_{0},\tau)>0 $ sufficiently small such that if
\begin{equation}
|\partial _{\lambda}^{l}P|_{D(r,s)\times \Lambda}\leq \epsilon
s^{2} \gamma ^{4m^{2}(n+1)},\ l \leq l_{0},
\end{equation}
then there exist Cantor sets $\Lambda_{\gamma}\subset \Lambda$
with $|\Lambda \backslash
\Lambda_{\gamma}|=O(\gamma^{\frac{1}{n-1}})$ and a $C^{l_{0}-1}$
Whitney smooth family of symplectic transformations
$$\Psi_{\lambda}:D(\frac{r}{2},\frac{s}{2})\longrightarrow D(r,s),
\ \lambda \in \Lambda_{\gamma}$$ such that
\begin{align}
H\circ
\Psi_{\lambda}(x,y,u)&=e_{\ast}(\lambda)+\langle\omega_{\ast}(\lambda),y\rangle
\nonumber\\
&+\frac{1}{2}\langle A_{\ast}(\lambda)y,y\rangle+\frac{1}{2}
\langle M_{\ast}(\lambda)u,u\rangle+P_{\ast}(x,y,u,\lambda)
\end{align}
and the following holds
$$(\omega_{\ast}(\lambda))_{i_{q}}\equiv
(\omega(\lambda))_{i_{q}},~q=1,\cdots, d.$$ Thus, all unperturbed
tori $T_{\lambda}=T^{n}\times \{0\}\times \{0\}$ with $\lambda \in
\Lambda_{\gamma}$ will persist and preserve the frequency
components $\omega_{i_{1}},\cdots,\omega_{i_{d}}$ of  the
unperturbed tangent frequencies $\omega(\lambda)$.

2) Assume $A(\lambda)$ is nonsingular on $\Lambda$ and let $\tau
>n-1$ be fixed. Then there exists an
$\epsilon=\epsilon(r,s,l_{0},\tau)>0$ sufficiently small such that
if (1.3) holds, then each torus $T_{\lambda}=T^{n}\times
\{0\}\times \{0\},\lambda \in \Lambda_{\gamma},$ will persist with
the normal form (1.4), and gives rise to an analytic, invariant
perturbed torus which preserves its tangent frequencies.
\end{theorem}

The $n_{0}$ in the above theorem can be arbitrary positive
integer. When $n_{0}\leq n$, Theorem 1.1 has applications to
Hamiltonian  system (1.1) with respect to the persistence of
invariant tori on sub-manifolds of $G$.

Consider (1.1) and let $S$ be an $n_{0}\ (\leq n)$ dimensional,
$C^{l_{0}}\ (l_{0}\geq n)$ sub-manifold of $G$ which is either
closed or with boundary. Denote
$$\omega(y)={\frac{\partial N}{\partial y}}(y),\ A(y)=
{\frac{\partial ^{2}N}{\partial y^{2}}}(y) ,\ y\in G.$$

We assume the following conditions:

\begin{itemize}
\item[{\bf A1)'}]
  For any
coordinate chart $(\phi,U)$ of $S$, rank$\{\frac {\partial
^{\alpha}(\omega\circ \phi^{-1})} {\partial
\lambda^{\alpha}}:|\alpha|\leq n-1\}=n$ for all $\lambda \in
\phi(U)\subset R^{n_{0}}.$

\item[{\bf A2)'}]  meas$\{\lambda: |i\langle
k,\omega\circ\phi^{-1}\rangle+\langle
l,\Omega\circ\phi^{-1}\rangle|=0\}=0, 0<|k|\leq K,|l|\leq 2,$ for
all $\lambda \in \phi(U)\subset R^{n_{0}}, $ where $K$ is defined
as in A2).

\item[{\bf A3)'}]  rank$A(y)\equiv d$  on $S$, and, there is a
smoothly varying, nonsingular, $d\times d$ principal minor
$\tilde{A}(y)$ of $A(y)$ on $S$.
\end{itemize}

\begin{corollary}Consider (1.1).

1) Assume A0),A1)', A2)', A3)' and let $\tau > n(n-1)-1$ be fixed.
Then there is an $\epsilon_{0}=\epsilon_{0}(r,s,l_{0},S,\tau)>0$
and a family of Cantor sets $S_{\epsilon}\subset S,0<\epsilon
\leq\epsilon _{0},$ with $|S\backslash
S_{\epsilon}|=O(\gamma^{\frac{1}{n-1}})$, such that for each
$y\in S_{\epsilon}$, the unperturbed torus $T_{y}$ persists and
gives rise to an analytic, invariant torus of the perturbed system
whose tangent frequencies $\omega_{\epsilon}$ satisfies
$$(\omega_{\epsilon}(y))_{i_{q}}=(\omega(y))_{i_{q}},q=1,\cdots,d,$$
where $i_{1},\cdots,i_{d}$ are the row indices (in the natural
order) of $\tilde {A}(y)$ located in $A(y)$. Moreover, these
perturbed tori form a Whitney smooth family.

2) Assume that $A(y)$ is nonsingular on $S$ and let $\tau >n-1$ be
fixed. Then each torus $T_{y},y\in S_{\epsilon}$, will persist and
gives rise to an analytic invariant perturbed torus with unchanged
tangent frequencies.

3) Let $y_{0}\in S_{\epsilon}$ in 1). Then (1.1) admits the
following normal form:
\begin{align*}
H_{y_{0}}(x,y,u)&=e_{\ast}(y_{0})+\langle\omega_{\ast}(y_{0}),y-y_{0}\rangle
\\
&+\frac{1}{2}\langle A_{\ast}(y_{0})(y-y_{0}),(y-y_{0})\rangle
+\frac{1}{2}\langle M_{\ast}(y_{0})u,u\rangle,
\end{align*}
where $\omega_{\ast}(y_{0})$ is the tangent frequencies of the
perturbed torus associated to $y_{0}$.
\end{corollary}

Similar to \cite{cl}, to generalize the standard isoenergetic KAM
theorem, we have to assume an additional sub-isoenergetic
non-degenerate condition besides the R\"{u}ssmann non-degeneracy
on an energy surface. More precisely, let $S$ be a sufficiently
smooth, relatively open, bounded subset of $\{N(y)=E \}$. We
assume A1)' on $S$ and also the following {\it sub-isoenergetic
non-degeneracy}:
\begin{itemize}
\item[{\bf A1)"}]
 {\it There is a smoothly varying
$d\times d$ principal minor $\tilde{A}(y)$ of $A(y)$ on $S$ such
that
$${\rm det}
     {\left(
     \begin{array}{cc}
        \tilde{A}(y) & \omega^{\ast}(y)\\
        \omega^{\ast}(y)^{\top}  & 0
     \end{array}
     \right)
      }\neq 0,
$$
where $\omega^{\ast}(y)=\frac{\partial N}{\partial y^{\ast}}(y),\
y^{\ast}= (y_{i_{1}},\cdots,y_{i_{d}})^\top,$ and
$i_{1},\cdots,i_{d}$ denote the row indices of $\tilde{A}(y)$ in
$A(y)$.}
\end{itemize}

\begin{theorem}
Consider (1.1). Let $S$ be a sufficiently smooth, relatively open,
bounded subset of $\{N(y)=E\}$.

1) Assume A0), A1)', A2)' on $S$ and let $\tau >n(n-1)-1$ be
fixed. Then there is an
$\epsilon_{0}=\epsilon_{0}(r,s,l_{0},m,S,\tau)>0$ and a family of
Cantor sets $S_{\epsilon}\subset S, 0<\epsilon\leq \epsilon_{0},$
with $|S\backslash S_{\epsilon}|=O(\gamma^{\frac{1}{n-1}})$, such
that for each $y\in S_{\epsilon}$, the unperturbed torus $T_{y}$
persists and gives rise to an analytic, invariant torus
$T_{\epsilon,y}$ of the perturbed system on the energy surface
$\{H(x,y,u)=E \}$. Moreover, these perturbed tori form a local
Whitney smooth family.

2) If A1)" also holds on $S$, then each perturbed torus
$T_{\epsilon,y}$ preserves the ratio of the $i_{1},\cdots,i_{d}$
components of its tangent frequencies $\omega_{\epsilon}$, i.e.,
$$[\omega_{\epsilon,i_{1}}:\cdots:\omega_{\epsilon,i_{d}}]
=[\omega_{i_{1}}:\cdots:\omega_{i_{d}}],$$ where
$\omega_{i_{j}},\omega_{\epsilon,i_{j}}$ are the $i_{j}$-th
components of unperturbed and perturbed tangent frequencies
respectively, for $j=1,2,\cdots,d.$

3) For $y_{0}\in S_{\epsilon}$, (1.1) admits the same normal form
as in part 3) of Corollary 1.1.
\end{theorem}

\noindent {\bf Remark 1.1.} 1) Our result is almost parallel to
Chow, Li and Yi \cite{cl}, i.e., we have the same results for
lower dimensional invariant tori as that of the standard
Hamiltonians  which have been shown in \cite{cl}.

2) In fact, in our case, we also can obtain arbitrarily prescribed
high ordered normal form of Hamiltonian similar to Chow, Li and Yi
\cite{cl}. To do so, we only need to change the iteration scheme a
little, but we do not do it for brief.

\section{KAM step}

\setcounter{equation}{0}
\renewcommand{\theequation}{\thesection.\arabic{equation}}

In this section, we describe the iterative scheme for the
Hamiltonian (1.2) in one KAM step. For simplicity, we set
$l_{0}=n.$

Consider (1.2) and initially set\\
$$r_{0}=r,\gamma_{0}=\gamma,s_{0}=s,\Lambda_{0}=
\Lambda,H_{0}=H,e_{0}=e,\omega_{0}=\omega,$$
$$A_{0}=A,\tilde{A}_{0}=\tilde{A},M_{0}=M,P_{0}=P,$$
$$N_{0}=e_{0}(\lambda)+\langle\omega_{0}(\lambda),y\rangle
+\frac{1}{2}\langle A_{0}(\lambda)y,y\rangle + \frac{1}{2}\langle
M_{0}(\lambda)u,u\rangle.$$ Without loss of generality, we assume
that $0<r_{0},\ s_{0},\ \gamma_{0}<1$  and $\tilde{A}_{0}$ is the
ordered $d\times d$  principal minor of $A_{0}$.

In what follows, the Hamiltonian without subscripts denotes the
Hamiltonian in $\nu$-th step, while those with subscripts ``+"
denotes the Hamiltonian of $(\nu+1)$-th step. And we shall use
``$<\cdot$" to denote ``$<c$" with a constant $c$ which is
independent of the iteration step. To simplify the notations, we
shall suspend the $\lambda$ dependence in most terms of this
section.

Suppose at the $\nu$-th step, we have arrived at the following
real analytic Hamiltonian:
\begin{align}
H&=N+P,
\nonumber\\
N&=e(\lambda)+\langle\omega(\lambda),y\rangle+\frac{1}{2} \langle
A(\lambda)y,y\rangle+\frac{1}{2}\langle M(\lambda)u,u\rangle,
\end{align}
which is defined on a phase domain $D(r,s)$ and depends smoothly
on $\lambda \in \Lambda$, where $\Lambda \subset \Lambda_{0}$.
Suppose that the $d \times d$ ordered principal minor $\tilde{A}$
of $A$ and $M$ are non-singular on $\Lambda$, and moreover,
$P=P(x,y,u,\lambda)$ satisfies
\begin{equation}
|\partial_{\lambda}^{l}P|_{D(r,s)}\leq  \epsilon
s^{2}\gamma^{4m^{2}(n+1)},|l|\leq n.
\end{equation}

We will construct a symplectic transformation $\Phi_{+}$, which
transforms the Hamiltonian (2.1) into the Hamiltonian of the next
KAM cycle (the $(\nu +1)$-th step), i.e.,
$$H_{+}=H\circ \Phi_{+}=N_{+}+P_{+},$$
where $N_{+}, P_{+}$ satisfy  similar conditions as $N, P$
respectively on $D(r_{+},s_{+})\times \Lambda_{+}.$

Define
\begin{align*}
\epsilon_{+}&=\epsilon^{\frac{10}{9}},
\\
r_{+}&=\frac{r}{2}+\frac{r_{0}}{4},
\\
s_{+}&=\frac{1}{8}\alpha s,\alpha =\epsilon^{\frac{1}{3}},
\\
\gamma_{+}&=\frac{\gamma}{2}+\frac{\gamma_{0}}{4},
\\
K_{+}&=([\log\frac{1}{\epsilon}]+1)^{a^{*}+2},
\\
D_{\frac{i}{8}\alpha}&=D(r_{+}+\frac{i-1}{8}(r-r_{+}),\frac{i}{8}\alpha
s),i=1,2,\cdots,8,
\\
D_{+}&=D_{\frac{1}{8}\alpha}=D(r_{+},s_{+}),
\\
\Lambda_{+}&=\{\lambda \in \Lambda:|i\langle
k,\omega(\lambda)\rangle+\langle
l,\Omega(\lambda)\rangle|>\frac{\gamma}{|k|^{\tau}},
\\
&|l|\leq 2,\ 0<|k|\leq K_{+}\},
\\
\Gamma(r-r_{+})&=\sum\limits_{0<|k|\leq
K_{+}}|k|^{\tau(n+1)4m^{2}+4m^{2}n}e^{-|k|\frac{r-r_{+}}{8}},
\end{align*}
where $a^{*}$ is a constant such that $(\frac{10}{9})^{a^{*}}>2$.

\subsection{Truncation}
Express $P$ into Taylor-Fourier series
$$P=\sum\limits_{k\in Z^{n},l\in Z_{+}^{n},p
\in Z_{+}^{2m}}p_{klp}y^{l}u^{p}e^{i<k,x>}$$ and let $R$ be the
truncation of $P$ with the form
\begin{equation}
R=\sum\limits_{|k|\leq K_{+},|l|+|p|\leq
2}p_{klp}y^{l}u^{p}e^{i<k,x>}.
\end{equation}

\begin{lemma} Assume that \\
{\bf H1)} $\int_{K_{+}}^{\infty}x^{n}e^{-x\frac{r-r_{+}}{8}}
{\rm d}x\leq \epsilon.$\\
Then the following
\begin{equation}
|\partial_{\lambda}^{l}(P-R)|_{D_{\frac{7}{8}\alpha}}\leq \cdot
\epsilon^{2}s^{2}\gamma^{4m^{2}(n+1)},
\end{equation}
$$|\partial_{\lambda}^{l}R|_{D_{\frac{7}{8}\alpha}}\leq
\cdot \epsilon s^{2}\gamma^{4m^{2}(n+1)}$$ hold for all $|l|\leq
n,\ \lambda \in \Lambda.$
\end{lemma}
{\bf Proof.} Let
$$I=\sum\limits_{|k|> K_{+},|l|+|p|\leq 2}p_{klp}y^{l}u^{p}e^{i<k,x>},$$
$$II=\sum\limits_{|k|\leq K_{+},|l|+|p|>2}p_{klp}y^{l}u^{p}e^{i<k,x>}.$$
Then
$$P-R=I+II.$$
By Cauchy estimate  and H1), we have
\begin{align*}
|\partial_{\lambda}^{l}I|_{D(r_{+}+\frac{7}{8}(r-r_{+}),s)}&\leq\sum_{|k|>K_{+}}|
\partial_{\lambda}^{l}P|_{D(r,s)}e^{-|k|\frac{r-r_{+}}{8}}
\\
&\leq \epsilon s^{2}\gamma^{4m^{2}(n+1)}
\sum_{x=K_{+}}^{\infty}x^{n}e^{-|k|\frac{r-r_{+}}{8}}
\\
&\leq \epsilon
s^{2}\gamma^{4m^{2}(n+1)}\int_{K_{+}}^{\infty}x^{n}e^{-x\frac{r-r_{+}}{8}}{\rm
d}x \leq \epsilon^{2}s^{2}\gamma^{4m^{2}(n+1)}.
\end{align*}
It follows that
\begin{align*}
|\partial_{\lambda}^{l}(P-I)|_{D(r_{+}+\frac{7}{8}(r-r_{+}),s)}&\leq
|\partial_{\lambda}^{l}P|_{D(r,s)}+|\partial_{\lambda}^{l}I|_{D(r_{+}
+\frac{7}{8}(r-r_{+}),s)}\\
&\leq \cdot \epsilon s^{2}\gamma^{4m^{2}(n+1)}.
\end{align*}

For $|q|=3$, let $\int$ be the obvious anti-derivative of
$\frac{\partial^{|l|+|p|}}{\partial y^{l}u^{p}}, |l|+|p|=3$. Then
by Cauchy estimate, it follows that
\begin{align*}
|\partial_{\lambda}^{l}II|_{D_{\frac{7}{8}\alpha}}&=|\partial_{\lambda}^{l}\int
\frac{\partial^{|l|+|p|}} {\partial
y^{l}u^{p}}\sum\limits_{|k|\leq K_{+},|l|+|p|>2}
p_{klp}y^{l}u^{p}e^{i<k,x>}{\rm d}y{\rm d}u|_{D_{\frac{7}{8}\alpha}}\\
&\leq|\frac{1}{s^{3}}\int |\partial_{\lambda}^{l}(P-I)|{\rm
d}y{\rm d}u
|_{D_{\frac{7}{8}\alpha}}\\
&\leq \cdot \frac{\alpha^{3}s^{3}}{s^{3}}|
\partial_{\lambda}^{l}(P-I)|_{D_{\frac{7}{8}\alpha}}\\
&\leq \cdot \epsilon^{2} s^{2}\gamma^{4m^{2}(n+1)}.
\end{align*}

Thus,
$$|\partial_{\lambda}^{l}(P-R)|_{D_{\frac{7}{8}\alpha}}
\leq \cdot \epsilon^{2}s^{2}\gamma^{4m^{2}(n+1)},$$ and therefore,
$$|\partial_{\lambda}^{l}R|_{D_{\frac{7}{8}\alpha}}
\leq |\partial_{\lambda}^{l}(P-R)|_{D_{\frac{7}{8}\alpha}}
+|\partial_{\lambda}^{l}P|_{D(r,s)}\leq \cdot \epsilon
s^{2}\gamma^{4m^{2}(n+1)}.\quad\Box$$

\subsection{Averaging and solving homogeneous equation}
To transform (2.1) into the Hamiltonian in the next KAM step, we
shall construct a  symplectic  transformation as the time 1-map
$\phi_{F}^{1}$ of the flow generated by a Hamiltonian $F$. To this
end, suppose $F$ has the following form:
\begin{align}
F&=\sum\limits_{0<|k|\leq K_{+},|l|+|p|\leq
2}f_{klp}y^{l}u^{p}e^{i\langle k,x\rangle}+\langle
f_{001},u\rangle+\langle f_{011}y,u\rangle
\nonumber \\
&=\sum\limits_{0<|k|\leq K_{+},|l|\leq 1}f_{klp}y^{l}e^{i\langle
k,x\rangle} +\sum\limits_{0<|k|\leq
K_{+},|l|=2}f_{klp}y^{l}e^{i\langle k,x\rangle}
\nonumber \\
&\quad+\sum\limits_{0\leq |k|\leq
K_{+},|l|=|p|=1}f_{klp}y^{l}u^{p}e^{i\langle k,x\rangle}
+\sum\limits_{0\leq |k|\leq K_{+},|p|=1}f_{klp}u^{p}e^{i\langle
k,x\rangle}
\nonumber \\
&\quad+\sum\limits_{0<|k|\leq K_{+},|p|=2}f_{klp}u^{p}e^{i\langle
k,x\rangle}
\nonumber \\
&\equiv F_{0}+F_{1}+F_{2}+F_{3}+F_{4},
\end{align}
which satisfies the equation
\begin{equation}
\{N,F\}+R-[R]+\langle p_{001},u\rangle+\langle
p_{011}y,u\rangle=0,
\end{equation}
where $[R]=\frac{1}{{(2\pi)}^{n}}\int_{T^{n}}R {\rm d}x.$
Substituting $N$ and (2.5) into (2.6) and comparing coefficients
yields that
\begin{align}
\triangle f_{kl0}&=-p_{kl0},0<|k|\leq K_{+},|l|\leq 1,\\
\triangle f_{k20}&=-p_{k20},0<|k|\leq K_{+},|l|=2,\\
(\triangle I_{2m}+MJ )f_{kl1}&=-p_{kl1},0\leq |k|\leq
K_{+},|l|=|p|=1,\\
(\triangle I_{2m}+MJ )f_{k01}&=-p_{k00,1},0 \leq |k|\leq
K_{+},|l|=0,|p|=1,
\\
(\triangle I_{2m}+MJ )f_{k02}-f_{k02}JM&=-p_{k02},0< |k|\leq
K_{+},|l|=0,|p|=2,
\end{align}
where $\triangle =i\langle k,\omega(\lambda)+A(\lambda)y\rangle,$
and $J$ is the standard symplectic matrix in $R^{2m}$.

\subsection{Estimate on $F$}
Let $\Omega_{j},\ j=1,\cdots,2m$ be  the eigenvalues of $JM$,
where $\Omega_{j}$ depends smoothly on $\lambda$. Then by the
non-degeneracy of $M$ there exists a constant $c$ such that
\begin{equation}
|\Omega_{j}|\geq c,\ j=1,\cdots,2m.
\end{equation}

\begin{lemma}Assume that\\
{\bf H2)} $2M_{\ast}s\leq \frac{\gamma}{K_{+}^{\tau+1}}$,\\
where $M_{\ast}$ is a constant defined to satisfy
$|A(\lambda)|\leq M_{\ast}$ on $\Lambda_{0}$. Then we have the
following result:
\begin{equation}
\frac{1}{s^{2}}|\partial_{\lambda}^{l}F|_{D(r_{+}+\frac{7}{8}(r-r_{+}),s)\times
\Lambda_{+}}\leq \cdot (\epsilon \Gamma+\epsilon).
\end{equation}
\end{lemma}
{\bf Proof.} By the definition of $\Lambda$ and H2) we obtain that
\begin{equation}
|\triangle|>\frac{\gamma}{2|k|^{\tau}}.
\end{equation}

And by the definition of $\Lambda_{+}$, H2), (2.12) and Lemma A.3,
we have the following results:

\begin{equation}
|det(\triangle I_{2m}-JM)|_{\Lambda_{+}}>\cdot
(\frac{\gamma}{|k|^{\tau}})^{2m}
\end{equation}

\begin{equation}
|det(\triangle I_{4m^{2}}-I_{2m}\otimes (JM)-(JM)\otimes
I_{2m})|_{\Lambda_{+}}
>\cdot (\frac{\gamma}{|k|^{\tau}})^{4m^{2}}.
\end{equation}

To estimate $F$, we must estimate
$\partial_{\lambda}^{l}\triangle^{-1}, \
\partial_{\lambda}^{l}(\triangle I_{2m}+MJ)^{-1}, \
\partial_{\lambda}^{l}(\triangle I_{4m^{2}}-I_{2m}\otimes
(JM)-(JM)\otimes I_{2m})^{-1}$ at first. For
$\partial_{\lambda}^{l}\triangle^{-1},$ we have
\begin{equation}
|\partial_{\lambda}^{l}\triangle^{-1}|_{\Lambda_{+}}\leq \cdot
|\triangle^{-1}|^{l+1} |k|^{l}\leq \cdot
\frac{|k|^{\tau(l+1)+l}}{\gamma^{l+1}},\ |l|\leq n.
\end{equation}
We note by the definition of $\Lambda_{+}$, H2), (2.12) and Lemma
A.3 that
\begin{equation}
|\partial_{\lambda}^{l}(\triangle
I_{2m}+MJ)^{-1}|_{\Lambda_{+}}\leq \cdot
(\frac{|k|^{\tau}}{\gamma})^{2m(l+1)}\times |k|^{2ml},
\end{equation}
\begin{equation}
|\partial_{\lambda}^{l}(\triangle I_{4m^{2}}-I_{2m}\otimes
(JM)-(JM)\otimes I_{2m})^{-1}|_{\Lambda_{+}} \leq \cdot
(\frac{|k|^{\tau}}{\gamma})^{4m^{2}(l+1)}\times |k|^{4m^{2}l}.
\end{equation}
Thus on $D(r,s)\times \Lambda_{+}$, we have
\begin{align}
|\partial_{\lambda}^{l}f_{kl0}|&=|\partial_{\lambda}^{l}(\triangle^{-1}p_{kl0})|
\leq \cdot \frac{|k|^{\tau(l+1)+l}}{\gamma^{l+1}}e^{-|k|r}\epsilon
s^{2-|l|}\gamma^{4m^{2}(n+1)},\\
|\partial_{\lambda}^{l}f_{k20}|&=|\partial_{\lambda}^{l}(\triangle^{-1}p_{k20})|
\leq \cdot \frac{|k|^{\tau(l+1)+l}}{\gamma^{l+1}}e^{-|k|r}\epsilon
\gamma^{4m^{2}(n+1)},
\\
|\partial_{\lambda}^{l}f_{k11}|&=|\partial_{\lambda}^{l}[(\triangle
I_{2m}+MJ)^{-1}p_{k11}]| \leq |\partial_{\lambda}^{l}[(\triangle
I_{2m}-JM)^{-1}p_{k11}]|
\nonumber\\
&\leq \cdot
\frac{|k|^{\tau(l+1)2m+2ml}}{\gamma^{(l+1)2m}}e^{-|k|r}\epsilon
\gamma^{4m^{2}(n+1)},
\\
|\partial_{\lambda}^{l}f_{k01}|&=|\partial_{\lambda}^{l}[(\triangle
I_{2m}+MJ)^{-1}p_{k01}]| \leq |\partial_{\lambda}^{l}[(\triangle
I_{2m}-JM)^{-1}p_{k01}]|
\nonumber\\
&\leq \cdot
\frac{|k|^{\tau(l+1)2m+2ml}}{\gamma^{(l+1)2m}}e^{-|k|r}\epsilon s
\gamma^{4m^{2}(n+1)},
\end{align}

\begin{align}
|\partial_{\lambda}^{l}f_{k02}|&=|\partial_{\lambda}^{l}[(\triangle
I_{4m^{2}}-I_{2m}\otimes (JM)-(JM)\otimes I_{2m})^{-1}p_{k02}]|
\nonumber \\
&\leq \cdot
\frac{|k|^{\tau(l+1)4m^{2}+4m^{2}l}}{\gamma^{(l+1)4m^{2}}}e^{-|k|r}\epsilon
\gamma^{4m^{2}(n+1)},
\end{align}
where $0<|k|\leq K_{+}.$ When $k=0$, we have
\begin{equation}
|f_{011}|\leq \cdot \epsilon \gamma^{4m^{2}(n+1)}, |f_{001}|\leq
\cdot \epsilon s\gamma^{4m^{2}(n+1)},
\end{equation}
on account of
\begin{equation}
MJf_{011}=-p_{011},\ MJf_{001}=-p_{001}.
\end{equation}
Therefore we obtain the estimate of $F$:
\begin{align}
\frac{1}{s^{2}}&|\partial_{\lambda}^{l}F|_{D(r_{+}+\frac{7}{8}(r-r_{+}),s)\times
\Lambda_{+}}
 \nonumber \\
&\leq\cdot\sum_{0<|k|\leq
K_{+}}\frac{|k|^{\tau(n+1)4m^{2}+4m^{2}n}}{\gamma^{(n+1)4m^{2}}}
e^{-|k|\frac{r-r_{+}}{8}}\epsilon \gamma^{4m^{2}(n+1)}
 \nonumber \\
&\quad+\cdot \epsilon \gamma^{4m^{2}(n+1)}
 \nonumber \\
&\leq \cdot \epsilon\Gamma+\cdot\epsilon.\quad\Box
\end{align}

Denote $D_{i}=D(r_{+}+\frac{i}{4}(r-r_{+}),\frac{i}{4}s),\
i=1,2,3.$ By (2.27) and Cauchy estimate, we obtain on $D_{3}
\times \Lambda_{+}$ that :
\begin{equation}
(r-r_{+})|\partial_{\lambda}^{l}F_{x}|,\
s|\partial_{\lambda}^{l}F_{y}|,\ s|\partial_{\lambda}^{l}F_{u}|
\leq \cdot \epsilon (\Gamma+1)s^{2}.
\end{equation}

Since $F$ is a polynomial of $y$ and $u$ with order 2, by (2.28)
we obtain
\begin{equation}
|D^{j}F|_{D_{2}\times \Lambda_{+}}\leq \cdot \epsilon(\Gamma+1),\
|j|\geq 2.
\end{equation}

Let $F$ be the Hamiltonian (2.5) with coefficients given by Lemma
2.2. Using $\phi_{F}^{t}$ denotes the flow generated by $F$, then
\begin{align}
H\circ \phi_{F}^{1}&=(N+R)\circ \phi_{F}^{1}+(P-R)\circ
\phi_{F}^{1}
\nonumber \\
&=N+\{N,F\}+R+\int_{0}^{1}\{R_{t},F \}\circ \phi_{F}^{t}{\rm
d}t+(P-R)\circ \phi_{F}^{1}
\nonumber \\
&=N+[R]-\langle p_{001},u\rangle-\langle p_{011}y,u\rangle
\nonumber \\
&\quad+\{N,F\}+R-[R]+\langle p_{001},u\rangle+\langle
p_{011}y,u\rangle
\nonumber \\
&\quad+\int_{0}^{1}\{R_{t},F \}\circ \phi_{F}^{t}{\rm
d}t+(P-R)\circ \phi_{F}^{1}
\nonumber \\
&=(N+[R]-\langle p_{001},u\rangle-\langle p_{011}y,u\rangle)
\nonumber \\
&\quad+(\int_{0}^{1}\{R_{t},F \}\circ \phi_{F}^{t}{\rm
d}t+(P-R)\circ \phi_{F}^{1})
\nonumber \\
&=\bar{N}_{+}+\bar{P}_{+},
\end{align}
where $R_{t}=(1-t)\{N,F\}+R$.

This completes the averaging process.

\subsection{Translation and partial non-degeneracy}

Denote by $Y,\ P_{010}$  the vectors formed by the first d
components of $y,\ p_{010}$ respectively. Then it is easy to see
that the equation
\begin{equation}
\tilde{A}Y=-P_{010}
\end{equation}
has a unique solution $Y^{\ast}$ on $D(r,s)$ which  depends
smoothly on $\lambda$. If we denote
$$y^{\ast}={
   \left(
   \begin{array}{c}
           Y^{\ast}\\
           0
   \end{array}
   \right)
   },$$
then by (2.31), we obtain
\begin{equation}
Ay^{\ast}=- {
   \left(
   \begin{array}{c}
           P_{010}\\
           0
   \end{array}
   \right)
   }.
\end{equation}

Consider the translation
$$\phi:x\rightarrow x,\ y\rightarrow y+y^{\ast},\ u\rightarrow u$$
and denote
$$\Phi_{+}=\phi_{F}^{1}\circ \phi.$$

Then
\begin{align}
H\circ \Phi_{+}&=N_{+}+P_{+},
\nonumber \\
N_{+}&=\bar{N}_{+}\circ \phi-\psi
\nonumber \\
&=e_{+}+\langle\omega_{+},y\rangle+\frac{1}{2}\langle
A_{+}y,y\rangle+\frac{1}{2}\langle M_{+}u,u\rangle,
\nonumber \\
P_{+}&=\bar{P}_{+}\circ \phi+\psi,
\end{align}
where
\begin{align}
e_{+}&=e+\langle\omega,y^{\ast}\rangle+\frac{1}{2}\langle
Ay^{\ast},y^{\ast}\rangle+[R](y^{\ast}),\\
\omega_{+}&=\omega+p_{010}- {
   \left(
   \begin{array}{c}
           P_{010}\\
           0
   \end{array}
   \right)
   }
,\\
M_{+}&=M+p_{002},\\
A_{+}&=A+\partial_{y}^{2}[R](y^{\ast}),\\
\psi&=\langle\partial_{y}[R](y^{\ast}),y\rangle-\langle
p_{010},y\rangle=2\langle p_{020}y^{\ast},y\rangle.
\end{align}

\subsection{Estimate on new normal form $N_{+}$}
\begin{lemma} We have the following holds for all $|l|\leq n$:
\begin{align}
&|\partial_{\lambda}^{l} y^{\ast}|_{\Lambda_{+}}
\leq \cdot \epsilon s \gamma^{4m^{2}(n+1)},\\
&|\partial_{\lambda}^{l} (e_{+}-e)|_{\Lambda_{+}}
\leq \cdot \epsilon s \gamma^{4m^{2}(n+1)},\\
&|\partial_{\lambda}^{l} (\omega_{+}-\omega)|_{\Lambda_{+}}
\leq \cdot \epsilon s \gamma^{4m^{2}(n+1)},\\
&|\partial_{\lambda}^{l} (A_{+}-A)|_{\Lambda_{+}}
\leq \cdot \epsilon \gamma^{4m^{2}(n+1)},\\
&|\partial_{\lambda}^{l} (M_{+}-M)|_{\Lambda_{+}}\leq \cdot
\epsilon \gamma^{4m^{2}(n+1)}.
\end{align}
\end{lemma}
{\bf Proof.} It is very clear by (2.31) and
(2.34)-(2.37).\quad$\Box$

\begin{lemma} Assume that \\
{\bf H3)} $\epsilon(\Gamma+1)s<\frac{1}{8}(r-r_{+});
\ \epsilon(\Gamma+1)s<\frac{1}{8}\alpha s.$\\
Then for all $0\leq t\leq 1$,
\begin{equation}
\Phi_{+}=\phi_{F}^{1}\circ
\phi:D_{+}=D_{\frac{1}{8}\alpha}\rightarrow
D_{\frac{1}{2}\alpha}\subset D(r,s),
\end{equation}
more precise,
\begin{align}
\phi :D_{\frac{1}{8}\alpha}\rightarrow D_{\frac{1}{4}\alpha},\\
\phi_{F}^{t}: D_{\frac{1}{4}\alpha}\rightarrow
D_{\frac{1}{2}\alpha}
\end{align}
are well defined, real analytic and depend smoothly on $\lambda
\in \Lambda_{+}$.
\end{lemma}
{\bf Proof.} (2.45) follows immediately from Lemma 2.3 and H3). To
prove (2.46), we rewrite
$\phi_{F}^{t}=(\phi_{F_{1}}^{t},\phi_{F_{2}}^{t},
\phi_{F_{3}}^{t})^{\top}$,where
$\phi_{F_{1}}^{t},\phi_{F_{2}}^{t}, \phi_{F_{3}}^{t}$ are
components of $\phi_{F}^{t}$ in the directions $x,y,u$
respectively. Let $(x,y,u)\in D_{\frac{1}{4}\alpha}$ and let
$t_{\ast}=Sup\{t\in [0,1]:\phi_{F}^{t}(x,y,u) \in
D_{\frac{6}{8}\alpha}\}$. Then for any $ 0 \leq t \leq t_{\ast}$,
\begin{align}
&|\phi_{F_{1}}^{t}(x,y,u)-x|\leq \int_{0}^{t}|F_{y}\circ
\phi_{F}^{s}|_{D_{\frac{6}{8}\alpha}}{\rm d}s \leq
|F_{y}|_{D_{\frac{6}{8}\alpha}}\leq \cdot
\epsilon(\Gamma+1)s<\frac{1}{8}(r-r_{+}) ,
 \nonumber \\
 &|\phi_{F_{2}}^{t}(x,y,u)-y|\leq \int_{0}^{t}|F_{x}\circ
 \phi_{F}^{s}|_{D_{\frac{6}{8}\alpha}}{\rm d}s
\leq |F_{x}|_{D_{\frac{6}{8}\alpha}}\leq \cdot
\epsilon(\Gamma+1)s^{2}<\frac{1}{8}\alpha s ,
 \nonumber \\
 &|\phi_{F_{3}}^{t}(x,y,u)-u|\leq \int_{0}^{t}|F_{u}\circ
 \phi_{F}^{s}|_{D_{\frac{6}{8}\alpha}}{\rm d}s
\leq |F_{u}|_{D_{\frac{6}{8}\alpha}}\leq \cdot
\epsilon(\Gamma+1)s<\frac{1}{8}\alpha s .
 \nonumber \\
\end{align}

It follows that $\phi_{F}^{t}(x,y,u)\in
D_{\frac{1}{2}\alpha}\subset D_{\alpha}.$ Thus, $t_{\ast}=1$ and
(2.46) holds.$\quad\Box$

Now we can give the estimate of $\Phi_{+}.$

\subsection{Estimate on the transformation $\Phi_{+}$}
\begin{lemma}For the transformation $\Phi_{+}$,
we have the following estimates:
\begin{equation}
|\Phi_{+}-id|_{D_{\frac{\alpha}{2}}}\leq \cdot
\epsilon(\Gamma+1)s, |D\Phi_{+}-Id|_{D_{\frac{\alpha}{2}}}\leq
\cdot \epsilon(\Gamma+1),
\end{equation}
where $id$ stands for the identity map, and $Id$ stands for the
elementary matrix.
\end{lemma}
{\bf Proof.} By
\begin{equation}
\phi_{F}^{1}=id+\int_{0}^{1}X_{F}\circ \phi_{F}^{s}{\rm d}s,
\end{equation}
we have
$$|\phi_{F}^{1}-id|\leq |X_{F}|_{D_{\frac{\alpha}{2}}}\leq
\cdot\epsilon(\Gamma+1)s.$$
For translation $\phi$ we have
\begin{equation}
|\phi-id|=|y^{\ast}|\leq \cdot \epsilon s\gamma^{4m^{2}(n+1)},
\end{equation}
so
$$|\phi|\leq 2.$$
Since
\begin{equation}
\Phi_{+}-id=(\phi_{F}^{1}-id)\circ \phi+ {
   \left(
   \begin{array}{c}
           0\\
           y^{\ast}\\
           0
   \end{array}
   \right)
   }
,
\end{equation}
we have
$$|\Phi_{+}-id|\leq \cdot\epsilon(\Gamma+1)s+\cdot \epsilon s\gamma^{4m^{2}(n+1)}
\leq \cdot\epsilon(\Gamma+1)s.$$ By (2.49) and (2.50), it follows
that
\begin{align*}
|D\phi_{F}^{1}-Id|&\leq 2|D^{2}F|\leq \cdot\epsilon(\Gamma+1),\\
|D\phi-Id|&\leq \cdot \epsilon\gamma^{4m^{2}(n+1)}.
\end{align*}
So by (2.51), we obtain the estimate of $D\Phi_{+}$:
\begin{align*}
&|D\Phi_{+}-Id|\\
&\leq |D(\phi_{F}^{1}-id)D\phi|+|Dy^{\ast}|\\
&\leq |D\phi_{F}^{1}-Id|\cdot|D\phi|+|Dy^{\ast}|\\
&\leq \cdot\epsilon(\Gamma+1).\quad\Box
\end{align*}

\subsection{Estimate on new perturbation $P_{+}$}
\begin{lemma}
Assume that\\
{\bf H4)} $\epsilon^{\frac{2}{9}}\Gamma \ll 1$,\\
then on $D_{+}\times \Lambda_{+}$,
\begin{equation}
|\partial_{\lambda}^{l}P_{+}|\leq
\epsilon_{+}s_{+}^{2}\gamma_{+}^{4m^{2}(n+1)},\ |l|\leq n.
\end{equation}
\end{lemma}
{\bf Proof.} Since
\begin{align*}
R_{t}=(1-t)\{N,F\}+R=tR+(1-t)[R]-(1-t)(<p_{001},u>+<p_{011}y,u>),
\end{align*}
it is easy to see that
$$|\partial_{\lambda}^{l} R_{t}|_{D(r,s)\times \Lambda_{+}}
\leq \cdot \epsilon s^{2}\gamma^{4m^{2}(n+1)}.$$ By the estimate
of $F$ and its derivative, we obtain that
\begin{align*}
|\partial_{\lambda}^{l} \{R_{t},F\}|_{D_{3}\times
\Lambda_{+}}&\leq |\partial_{\lambda}^{l}R_{tx}F_{y}|
+|\partial_{\lambda}^{l}R_{ty}F_{x}|+|\partial_{\lambda}^{l}R_{tu}F_{u}|\\
&\leq \cdot \epsilon^{2} s^{2} (\Gamma+1)\gamma^{4m^{2}(n+1)}.
\end{align*}
By Lemma 2.1, (2.38) and (2.39), we have on
$D_{\frac{\alpha}{2}}\times \Lambda_{+}$ the following holds
\begin{align*}
|\partial_{\lambda}^{l}(P-R)\circ \phi_{F}^{1}|&\leq \cdot
\epsilon^{2} s^{2}
 \gamma^{4m^{2}(n+1)}\\
 |\partial_{\lambda}^{l}\phi|&\leq \cdot \epsilon s\gamma^{4m^{2}(n+1)}\\
 |\partial_{\lambda}^{l}\psi|&\leq \cdot \epsilon^{2} s^{2}\gamma^{8m^{2}(n+1)}.
\end{align*}
So we have
$$|\partial_{\lambda}^{l}P_{+}|_{D_{+}\times \Lambda_{+}}
\leq \cdot \epsilon^{2}s^{2}\gamma^{4m^{2}(n+1)}(\Gamma +3)$$ by
the above estimate and (2.33). Thus it is enough to verify
\begin{equation}
\epsilon^{2}s^{2}\gamma^{4m^{2}(n+1)}(\Gamma +3)\leq
\epsilon_{+}s_{+}^{2}\gamma_{+}^{4m^{2}(n+1)}.
\end{equation}
By the definition of $\epsilon_{+},s_{+},\gamma_{+}$ and H4), it
is clear that it does hold.\quad$\Box$

This completes one cycle of KAM steps.

\section{Proof of main results}

\setcounter{equation}{0}
\renewcommand{\theequation}{\thesection.\arabic{equation}}

\subsection{Iteration lemma}
Considering (1.2), we define the following sequences inductively
for all $\nu=1,2,\cdots:$
\begin{align*}
r_{\nu}&=\frac{r_{\nu-1}}{2}+\frac{r_{0}}{4},
\\
s_{\nu}&=\frac{1}{8}\alpha_{\nu-1} s_{\nu-1},\alpha_{\nu}
=\epsilon_{\nu}^{\frac{1}{3}},
\\
\gamma_{\nu}&=\frac{\gamma_{\nu-1}}{2}+\frac{\gamma_{0}}{4},
\\
\epsilon_{\nu}&=\epsilon_{\nu-1}^{\frac{10}{9}},
\\
K_{\nu}&=([\log\frac{1}{\epsilon_{\nu-1}}]+1)^{a^{*}+2},
\\
D_{\frac{i}{8}\alpha}&=D(r_{+}+\frac{i-1}{8}(r-r_{+}),\frac{i}{8}\alpha
s),i=1,2,\cdots,8,
\\
D_{\nu}&=D(r_{\nu},s_{\nu}),
\\
\Gamma_{\nu}&=\Gamma(r_{\nu}-r_{\nu+1}),
\\
\Lambda_{\nu}&=\{\lambda \in \Lambda_{\nu-1}:|i\langle
k,\omega_{\nu-1}(\lambda)\rangle+\langle l,\Omega_{\nu
-1}(\lambda)\rangle|
\\
&\quad>\frac{\gamma_{\nu-1}}{|k|^{\tau}},\ |l|\leq 2,\ 0<|k|\leq
K_{\nu} \},
\end{align*}
where $a^{*}$ is a constant such that $(\frac{10}{9})^{a^{*}}>2.$

\begin{lemma}
If (1.3) holds for a sufficiently small $\epsilon_{0}$, then the
following holds for all $|l|\le n; \nu=1,2,\cdots.$

1)\begin{align} &|\partial_{\lambda}^{l} (e_{\nu}-e_{\nu
-1})|_{\Lambda_{\nu}}\leq \cdot \epsilon_{\nu -1} s_{\nu -1}
\gamma_{\nu -1}^{4m^{2}(n+1)},\\
&|\partial_{\lambda}^{l} (e_{\nu}-e_{0})|_{\Lambda_{\nu}}\leq
\cdot \epsilon_{0} s_{0}
\gamma_{0}^{4m^{2}(n+1)},\\
&|\partial_{\lambda}^{l} (\omega_{\nu}-\omega_{\nu
-1})|_{\Lambda_{\nu}}\leq \cdot \epsilon_{\nu -1}
s_{\nu -1} \gamma_{\nu -1}^{4m^{2}(n+1)},\\
&|\partial_{\lambda}^{l}
(\omega_{\nu}-\omega_{0})|_{\Lambda_{\nu}}\leq \cdot \epsilon_{0}
s_{0} \gamma_{0}^{4m^{2}(n+1)},\\
&|\partial_{\lambda}^{l} (A_{\nu}-A_{\nu -1})|_{\Lambda_{\nu}}\leq
\cdot \epsilon_{\nu -1}
\gamma_{\nu -1}^{4m^{2}(n+1)},\\
&|\partial_{\lambda}^{l} (A_{\nu}-A_{0})|_{\Lambda_{\nu}}\leq
\cdot \epsilon_{0}
\gamma_{0}^{4m^{2}(n+1)},\\
&|\partial_{\lambda}^{l} (M_{\nu}-M_{\nu -1})|_{\Lambda_{\nu}}\leq
\cdot \epsilon_{\nu -1}
\gamma_{\nu -1}^{4m^{2}(n+1)},\\
&|\partial_{\lambda}^{l} (M_{\nu}-M_{0})|_{\Lambda_{\nu}}\leq
\cdot \epsilon_{0}
\gamma_{0}^{4m^{2}(n+1)},\\
&|\partial_{\lambda}^{l} P_{\nu}|_{D_{\nu}\times
\Lambda_{\nu}}\leq \epsilon_{\nu}
s_{\nu}^{2}\gamma_{\nu}^{4m^{2}(n+1)}.
\end{align}

2) $(\omega_{\nu}(\lambda))_{q}=(\omega_{\nu -1}(\lambda))_{q}$
for all $q=1,2,\cdots d$ and $\lambda \in \Lambda_{\nu}.$

3) $\Phi_{\nu}:D_{\nu}\times \Lambda_{\nu}\rightarrow D_{\nu-1}$
is symplectic for each $\lambda \in \Lambda$, and
$$|\Phi_{\nu}-id|_{D_{\nu}\times \Lambda_{\nu}}
\leq\cdot\epsilon_{\nu-1}(\Gamma_{\nu-1}+1)s_{\nu-1}.$$ Moreover,
on $D_{\nu}\times \Lambda_{\nu},$
$$H_{\nu}=H_{\nu-1}\circ \Phi_{\nu}=N_{\nu}+P_{\nu},$$
where
\begin{align*}
H_{\nu}&=N_{\nu}+P_{\nu},\\
N_{\nu}&=e_{\nu}+\langle \omega_{\nu},y\rangle+\frac{1}{2}\langle
A_{\nu}y,y\rangle +\frac{1}{2}\langle M_{\nu}u,u\rangle,
\end{align*}
$A_{\nu}$ is real symmetric with its $d\times d$ ordered principal
minor $\tilde{A}_{\nu}$ being nonsingular on $\Lambda_{\nu}$.
\end{lemma}
{\bf Proof.} We only have to verify H1)-H4) for all $\nu$. For
simplicity, we let $r_{0}=1$.

First, we verify H1). By the choice of
$K_{+}=([log\frac{1}{\epsilon}]+1)^{a^{\ast}+2},$ where $a^{\ast}$
is a constant such that $(\frac{10}{9})^{a^{\ast}}>2$, we have
\begin{align*}
&\log(n+1)!+n(a^{\ast}+2)\log([\log\frac{1}{\epsilon}]+1)
-\frac{1}{2^{\nu+5}}(\log\frac{1}{\epsilon})^{a^{\ast}+2}\\
&\leq^{(I^{*})} \log(n+1)!+n(a^{\ast}+2)\log(\log\frac{1}
{\epsilon}+2)-(\log\frac{1}{\epsilon})^{2}\\
&\leq -\log\frac{1}{\epsilon}
\end{align*}
if $\epsilon_{0}$ is sufficiently small, where $(I^{*})$ holds
since $\frac{1}{2^{\nu+5}}(\log\frac{1}{\epsilon})^{a^{\ast}}\geq
1$ by the choice of $a^{\ast}$. Thus H1) holds since
$\int_{K_{+}}^{\infty}x^{n}e^{-x\frac{r-r_{+}}{8}}{\rm d}x\leq
(n+1)! K_{+}^{n}e^{-\frac{K_{+}}{2^{\nu+5}}r_{0}}$.

Then we verify H2). We have
\begin{align*}
2M_{\ast}sK_{+}^{\tau+1}&=2M_{\ast}(\frac{1}{8})^{\nu}s_{0}
\epsilon_{0}^{3[(\frac{10}{9})^{\nu}-1]}
(\log\frac{1}{\epsilon})^{(a^{\ast}+2)(\tau+1)}\\
&=2M_{\ast}(\frac{10}{9})^{\nu(a^{\ast+2})(\tau+1)}(\frac{1}{8})^{\nu}s_{0}
\epsilon_{0}^{3[(\frac{10}{9})^{\nu}-1]}
(\log\frac{1}{\epsilon_{0}})^{(a^{\ast}+2)(\tau+1)}\\
&=2M_{\ast}[(\frac{10}{9})^{(a^{\ast}+2)(\tau+1)}\frac{1}{8}]^{\nu}s_{0}
\epsilon_{0}^{3[(\frac{10}{9})^{\nu}-1]}
(\log\frac{1}{\epsilon_{0}})^{(a^{\ast}+2)(\tau+1)}\\
&\leq^{(II^{*})}C^{\nu}2M_{\ast}s_{0}\epsilon_{0}^{2[(\frac{10}{9})^{\nu}-1]}\\
&\leq^{(III^{*})} s_{0}\epsilon_{0}^{[(\frac{10}{9})^{\nu}-1]}\\
&\leq^{(IV^{*})}\frac{\gamma_{0}}{2},
\end{align*}
where $C=(\frac{10}{9})^{(a^{\ast}+2)(\tau+1)}\frac{1}{8}$,
$(II^{*}),\ (III^{*})$ can hold if $\epsilon_{0}$ is chosen
sufficiently small such that
$\epsilon_{0}^{[(\frac{10}{9})^{\nu}-1]}(\log\frac{1}{\epsilon_{0}})
^{(a^{\ast}+2)(\tau+1)}\leq 1$ and
$2M_{\ast}C^{\nu}\epsilon_{0}^{(\frac{10}{9})^{\nu}-1}\leq 1$, and
$(IV^{*})$ is easily done, say, set $s_{0}=\frac{\gamma_{0}}{2}$.

To verify verify H4), we note that
\begin{align*}
\epsilon^{\frac{2}{9}}\Gamma &\leq\epsilon^{\frac{2}{9}}
\int_{1}^{\infty}\lambda^{\tau(n+1)4m^{2}+4m^{2}n+n}
e^{-\lambda\frac{1}{2^{\nu +5}}}{\rm d} \lambda\\
&\leq \epsilon^{\frac{2}{9}}[\tau(n+1)4m^{2}+4m^{2}n+n+1]{\rm !}
2^{(\nu +5)[\tau(n+1)4m^{2}+4m^{2}n+n+1]}\\
&\leq \cdot \epsilon^{\frac{2}{9}}2^{\nu[\tau(n+1)4m^{2}+4m^{2}n+n+1]}\\
&\leq \cdot \epsilon_{0}^{\frac{2}{9}(\frac{10}{9})^{\nu}}
2^{\nu[\tau(n+1)4m^{2}+4m^{2}n+n+1]}\\
&\ll^{(V^{*})}1,
\end{align*}
where $(V^{*})$ holds if $\epsilon_{0}$ is chosen sufficiently
small.

Now, the rest work is to prove H3). By H4), we have
$$\epsilon(\Gamma+1)s\leq \epsilon^{\frac{7}{9}}s\leq \epsilon_{0}
^{\frac{7}{9}(\frac{10}{9})^{\nu}}s_{0}\leq^{(VI^{*})}
\frac{1}{2^{\nu+5}} =\frac{1}{8}(r-r_{+}),$$ where $(VI^{*})$
holds if $\epsilon_{0}$ is chosen sufficiently small. It is very
clear that $\epsilon(\Gamma+1)s<\frac{1}{8}\alpha s$ if H4) holds.

In the process of proof of the lemma, we have used the sufficient
smallness of $\epsilon_{0}$ in $(I^{*})-(VI^{*})$. In fact, the
existence of $\epsilon_{0}$ is obvious in  $(I^{*})-(VI^{*})$ in
spite that we do not give the explicit form. So, in the end, we
can choose the smallest $\epsilon_{0}$  of $(I^{*})-(VI^{*})$ as
the $\epsilon_{0}$ we need. This completes the proof of the
lemma.\quad$\Box$

\subsection{Proof of Theorem 1.1}
Denote
$$\Psi^{\nu}=\Phi_{1}\circ\Phi_{2}\circ\cdots\circ\Phi_{\nu},\nu=1,2,\cdots.$$
Then $\Psi^{\nu}:D_{\nu}\times \Lambda_{\nu}\rightarrow D_{0},$
and,
$$H_{0}\circ\Psi^{\nu}=H_{\nu}=N_{\nu}+P_{\nu},\nu=1,2,\cdots,$$
where $\Psi^{0}=id$. Let
$$\Lambda_{\ast}=\bigcap\limits_{\nu\geq 0}\Lambda_{\nu}.$$
Then $\Lambda_{\ast}$ is a Cantor-like set. First, we show that we
have the estimate
$$|\Lambda_{0}\backslash\Lambda_{\ast}|=O(\gamma^{\frac{1}{n-1}}),$$
we will divide the proof of which into two cases.\\
{\bf Case 1:}$n_{0}=n$.\\
According to \cite{xy},
$\{\partial^{\beta}\omega/\partial\lambda^{\beta}:\ \forall\beta,\
|\beta|=r\}$ and $\{D_{V}^{r}\omega:\forall V\in R^{n}\}$ are
linearly equivalent, where $r>0$ is an integer and
$D_{V}^{r}\omega=d^{r}/dt^{r}\omega(\lambda+tV)|_{t=0}.$ Since
(3.4) is satisfied by the extended tangent frequencies
$\omega_{\nu}$ on $\Lambda_{0},$
 A1) implies that if $\epsilon_{0}$ is sufficiently small, then
$${\rm rank}\{{{\partial^{\alpha}\omega_{\nu}}\over {\partial \lambda
^{\alpha}}}:|\alpha |\leq n-1\}=n$$ for all
$\lambda\in\Lambda_{0},\ \nu=0,1,\cdots.$ In the following proof,
we will omit the subscript $\nu.$ So there exist n integers $0\leq
r_{1},\cdots,r_{n}\leq n-1$ and n vectors $V_{1},\cdots,V_{n}\in
R^{n}$ such that
$${\rm rank}\{D_{V_{1}}^{r_{1}}\omega,\cdots,D_{V_{n}}^{r_{n}}\omega \}=n,
\ \forall \lambda \in \Lambda.$$ Denote
$B=(D_{V_{1}}^{r_{1}}\omega,\cdots,D_{V_{n}}^{r_{n}}\omega).$ Then
there exist a constant $\sigma>0$ such that for all
$(\lambda,V)\in\Lambda\times U,$
$$|BV|\geq \sigma,$$
where $U=\{V\in R^{n}:|V_{1}|+\cdots+|V_{n}|=1\}.$ Then it follows
that for some $1\leq i\leq n$ and $\forall k\in
Z^{n}\backslash\{0\},$
$$|\langle D_{V_{i}}^{r_{i}}\omega,\frac{k}{|k|}\rangle|\geq \frac{\sigma}{n}.$$
So by the definition of $K$, when $|k|>K$, we have
\begin{align*}
&|D_{V_{i}}^{r_{i}}(i\langle\frac{k}{|k|},\omega\rangle+
\frac{1}{|k|}\langle l,\Omega\rangle)|\\
&\geq \frac{\sigma}{n}-\frac{1}{|k|}|D_{V_{i}}^{r_{i}}\langle l,\Omega\rangle|\\
&\geq \frac{\sigma}{n}-\frac{2}{|k|}\frac{K}{\frac{4n}{\sigma}}\\
&\geq \frac{\sigma}{2n}.
\end{align*}
Let
\begin{align*}
R_{kV_{i}}&=\{t:|i\langle\frac{k}{|k|},\omega(\lambda+tV_{i})\rangle
+\frac{1}{|k|}\langle l,\Omega(\lambda+tV_{i})\rangle|\\
&\leq \frac{\gamma}{|k|^{\tau+1}},\lambda\in\Lambda,
\lambda+tV_{i}\in\Lambda \},\\
R_{k,l}&=\{\lambda\in\Lambda:|i\langle\frac{k}{|k|},
\omega(\lambda)\rangle+\frac{1}{|k|}\langle
l,\Omega(\lambda)\rangle|
 \leq \frac{\gamma}{|k|^{\tau+1}}\},\ |l|\leq 2.
\end{align*}
By Lemma A.1, when $|k|>K$, we have
$$|R_{kV_{i}}|\leq\cdot(\frac{\gamma}{|k|^{\tau+1}})^{\frac{1}{r_{i}}}
\leq\cdot(\frac{\gamma}{|k|^{\tau+1}})^{\frac{1}{n-1}}.$$ Then it
follows that
$$|R_{k,l}|\leq\cdot({\rm diam} ~\Lambda)^{n-1}\frac{\gamma^{\frac{1}{n-1}}}
{|k|^{\frac{\tau+1}{n-1}}}.$$ When $|k|\leq K$, by A2) we have
that $|R_{k,l}|\rightarrow 0(\gamma \rightarrow 0)$, i.e.,
$|R_{k,l}|=O(\gamma^{\frac{1}{n-1}}),(\gamma \rightarrow 0).$ So
we obtain that
\begin{align*}
|\Lambda_{0}\backslash \Lambda_{\ast}|&=
|\bigcup\limits_{k,l}R_{k,l}|\leq\sum\limits_{k,l}|R_{k,l}|\\
&\leq\cdot\gamma^{\frac{1}{n-1}}\sum\limits_{|k|>K}\frac{1}
{|k|^{\frac{\tau+1}{n-1}}}
+O(\gamma^{\frac{1}{n-1}})\sum\limits_{0<l\leq K}l^{n}\\
&=O(\gamma^{\frac{1}{n-1}}),
\end{align*}
which is the result we desired.\\
{\bf Case 2:} $n_{0}<n$. Let $\bar{\Lambda}=[1,2]^{n-n_{0}}$ and
define
\begin{align*}
&\tilde{\Lambda}=\Lambda_{0}\times \bar{\Lambda},\\
&\tilde{\Lambda}_{\ast}=\Lambda_{\ast}\times \bar{\Lambda},\\
&\tilde{\lambda}=(\lambda,\bar{\lambda})^{\top},\bar{\lambda}\in \bar{\Lambda},\\
&\tilde{\omega}_{\nu}(\tilde{\lambda})=\omega_{\nu}(\lambda),\nu=0,1,\cdots;\
\tilde{\lambda}\in \tilde{\Lambda}.
\end{align*}
Then by A1) it is clear that
$${\rm rank}\{\frac{\partial^{\alpha}\tilde{\omega}_{\nu}}
{\partial\tilde{\lambda}^{\alpha}}:\alpha\leq n-1\}=n$$ on
$\tilde{\Lambda}$ for all $\nu=0,1,\cdots,$ as $\epsilon_{0}$ is
sufficiently small. Similar to Case 1, we have that
$$|\tilde{\Lambda}\backslash\tilde{\Lambda}_{\ast}|=O(\gamma^{\frac{1}{n-1}}).$$

By Fubini's theorem,
$$|\Lambda_{0}\backslash\Lambda_{\ast}|=O(\gamma^{\frac{1}{n-1}})$$
as desired.

Since we mainly care about the persistence of invariant tori on
sub-manifolds, the measure estimate's case when $n_{0}>n$ is
omitted. In fact the reader can also see the reference Chow, Li
and Yi \cite{cl} or Li and Yi \cite{ly} for details.

Then we show the convergence of $H_{\nu}$ and $\Psi^{\nu}$.
Similar to the argument in \cite{cl} and \cite{ly}, in view of
Lemma 2.5 and Lemma 3.1, it concludes that $\Psi^{\nu}$ converges
uniformly to $\Psi^{\infty}$ on $D_{\infty}\times\Lambda_{\ast}$,
and under the map $\Psi^{\infty}$, $N_{\nu}$ converges uniformly
to $N_{\infty}$ on $D_{\infty}\times\Lambda_{\ast}$ with
$$N_{\infty}=e_{\infty}(\lambda)+\langle\omega_{\infty}(\lambda),y\rangle
+\frac{1}{2}\langle
A_{\infty}(\lambda)y,y\rangle+\frac{1}{2}\langle
M_{\infty}(\lambda)u,u\rangle.$$ Hence for each
$\lambda\in\Lambda_{\ast},\ T^{n}\times\{0\}\times\{0\}$ is an
analytic invariant torus of $H_{\infty}$ with the frequencies
$\omega_{\infty}(\lambda),$ which, by Lemma 3.1 2), satisfies
$$(\omega_{\infty}(\lambda))_{q}\equiv(\omega_{0}(\lambda))_{q},1\leq q\leq d.$$

Denote $\Psi_{\lambda}=\Psi^{\infty}(\cdot,\lambda)$, then
$\{\Psi_{\lambda}:\lambda \in \Lambda_{\ast} \}$ is a $C^{n-1}$
Whitney smooth family of analytic symplectic transformations on
$D(\frac{r_{0}}{2},\frac{s_{0}}{2})$ (see \cite{cl} for details).

Similar to \cite{cl}, following the Whitney extension of
$\Psi^{\nu}$'s, all $e_{\nu},\omega_{\nu},
A_{\nu},M_{\nu},P_{\nu},\nu=0,1,\cdots,$ admit uniform
$C^{n-1+\sigma_{0}}$ extensions in $\lambda\in\Lambda_{0}$ with
derivatives in $\lambda$ up to order $n-1$ satisfying the same
estimates (3.1)-(3.9). Thus, $e_{\infty},\omega_{\infty},
A_{\infty},M_{\infty},P_{\infty},$ are $C^{n-1}$ Whitney smooth in
$\lambda\in\Lambda_{\ast}$, and, the derivatives of
$(e_{\infty}-e_{0}),
(\omega_{\infty}-\omega_{0}),(A_{\infty}-A_{0}),(M_{\infty}-M_{0})$
satisfy similar estimates as in (3.2), (3.4), (3.6), (3.8).
Henceforth, the perturbed tori form a $C^{n-1}$ Whitney smooth
family on $\Lambda_{\ast}$.

This completes the proof of Theorem 1.1.\quad$\Box$

\subsection{Proof of Corollary 1.1}
{\bf Proof.} Without loss of generality, we assume that $S$ admits
a global coordinate, i.e., there is a bounded closed region
$\Lambda\subset R^{n_{0}}$ and a $C^{l_{0}}$ diffeomorphism $y$:
$\Lambda\rightarrow S$ such that $S=y(\Lambda)$. Let
$\lambda\in\Lambda$ and consider the transformation
$$y\rightarrow y+y(\lambda).$$
Then (1.1) gives rise to
\begin{align*}
H(x,y,u,\lambda)&=e(\lambda)+\langle\omega(\lambda),y\rangle
\\
&+\frac{1}{2}\langle A(\lambda)y,y\rangle +\frac{1}{2}\langle
M(\lambda)u,u\rangle+P(x,y,u,\lambda),
\end{align*}
where
\begin{align*}
e(\lambda)&=N(y(\lambda),0),\\
\omega(\lambda)&=\frac{\partial N}{\partial y}(y(\lambda),0),\\
A(\lambda)&=\frac{\partial^{2}N}{\partial y^{2}}(y(\lambda),0),\\
M(\lambda)&=\frac{\partial^{2}N}{\partial u^{2}}(y(\lambda),0),\\
P(x,y,u,\lambda)&=P(x,y+y(\lambda),u)+O(|y|^{3}+|u|^{3}).
\end{align*}
By the analysis of the Hamiltonian and assumption A0), there is no
$O(|yu|)$ in new perturbation $P$.

Let $s_{0}=\epsilon_{0}\gamma_{0}^{4m^{2}(n+1)}$. Then (1.3) holds
and the Corollary follows immediately from the Theorem 1.1 as
$\epsilon_{0}$ is sufficiently small.\quad$\Box$

\subsection{Proof of Theorem 1.2}
{\bf Proof.} By choosing $\lambda,\ \Lambda$ as in the Section 3.3
with the present $S$, the proof of Theorem 1.2 essentially follows
that of Theorem 1.1, except the translation
$$\phi:x\rightarrow x,\ y\rightarrow y+y^{\ast},\ u\rightarrow u$$
in Section 2.4 should be defined for purpose of eliminating the
energy drift at each KAM step. The rest proof is similar to
\cite{cl}.\quad$\Box$

\section{Some Examples}
In this section we give some examples to illustrate our results.\\
{\large\bf Example 4.1.}  We consider the following unperturbed
system:
$$N(y,u)=y_{1}+\frac{1}{2}y_{2}^{2}+\frac{\sqrt[]{2}}{2}(u^{2}+v^{2}),$$
where $y_{1},y_{2},u,v\in R^{1},$ that is $n=2,m=1.$ It is easy to
see that :
\begin{align*}
\Omega={ \left(
\begin{array}{c}
   i\sqrt[]{2}\\
   -i\sqrt[]{2}
\end{array}
\right) }, \omega={ \left(
\begin{array}{c}
   1 \\
   y_{2}
\end{array}
\right) }, A={ \left(
\begin{array}{cc}
   0 & 0\\
   0 & 1
\end{array}
\right) },
\\
R={ \left(
\begin{array}{cc}
   1 & 0\\
   y_{2} & \partial_{\lambda}y_{2}
\end{array}
\right) }, M={ \left(
\begin{array}{cc}
   \sqrt[]{2} & 0\\
   0 & \sqrt[]{2}
\end{array}
\right) },
\end{align*}
where $R$ stands for the matrix
$\{\frac{\partial^{\alpha}\omega}{\partial\lambda^{\alpha}}:
\alpha\leq n-1\}.$

1) We consider the persistence of invariant tori on the line
segment:
$$S_{1}:\ y_{1}(\lambda)=a_{1}\lambda,\ y_{2}(\lambda)=a_{2}
\lambda,\ \lambda\in[1,2].$$
Obviously A1)' holds on $S_{1}$ if and only if $a_{2}\neq 0$. We
can easily verify that
$$\{\lambda:|i\langle k,\omega\rangle+\langle l,\Omega\rangle|=
0,\ 0<|k|\leq K,\ |l|\leq 2\}$$
contains at most an isolated point, that is, A2)' holds. So by the
expression of $A$ and our Corollary 1.1 1), the majority 2-tori on
$S_{1}$ will persist with unchanged second component of tangent
frequencies. Since $A$ is  singular, part 2) of Corollary 1.1 is
not applicable.

2) We consider the persistence of invariant tori on the parabola:
$$S_{2}:\ y_{1}(\lambda)=a_{1}\lambda,\ y_{2}(\lambda)=a_{2}
\lambda^{2},\ \lambda\in[1,2].$$
Similar to 1), we can verify that A1)' holds if and only if
$a_{2}\neq 0$. And similar to 1), we can verify that
$$\{\lambda:|i\langle k,\omega\rangle+\langle l,\Omega\rangle|
=0,\ 0<|k|\leq K,\ |l|\leq 2\}$$
contains at most two points, i.e., A2)' holds. Also we obtain that
the majority 2-tori on $S_{2}$ will persist with unchanged second
component of  tangent frequencies. As $A$ is  singular, part 2) of
the Corollary 1.1 is not applicable.

Since the eigenvalues of $JM$ are pure imaginary, the persistent
tori are elliptic. \\
{\large\bf Example 4.2.}   We consider the following unperturbed
system:
$$N(y,u)=\frac{1}{2}y_{1}^{2}+\frac{1}{3}y_{2}^{3}+
\frac{\sqrt[]{2}}{2}(u_{1}^{2}+v_{1}^{2})
+\frac{\sqrt[]{3}}{2}(u_{2}^{2}+v_{2}^{2})$$ with
\begin{align*}
\omega={ \left(
\begin{array}{c}
   y_{1}\\
   y_{2}^{2}
\end{array}
\right) }, A={ \left(
\begin{array}{cc}
   1 & 0\\
   0 & 2y_{2}
\end{array}
\right) }, R={ \left(
\begin{array}{cc}
   y_{1} & \partial_{\lambda}y_{1}\\
   y_{2}^{2} & 2y_{2}\partial_{\lambda}y_{2}
\end{array}
\right) },
\\
\Omega={ \left(
\begin{array}{c}
   i\sqrt[]{2}\\
   -i\sqrt[]{2}\\
   i\sqrt[]{3}\\
   -i\sqrt[]{3}
\end{array}
\right) }, M={ \left(
\begin{array}{cccc}
   \sqrt[]{2} & 0 & 0 & 0\\
   0& \sqrt[]{2} &0&0\\
   0& 0& \sqrt[]{3}&0\\
  0& 0&0 & \sqrt[]{3}
\end{array}
\right) }.
\end{align*}

1)  We consider the persistence of invariant tori on the line
segment:
$$S_{1}:\ y_{1}(\lambda)=a_{1}\lambda,\ y_{2}(\lambda)=a_{2}\lambda,
\ \lambda\in[1,2].$$
It is easy to see that $R$ is nonsingular on $S_{1}$, i.e. A1)'
holds, if and only if $a_{1}a_{2}\neq 0$. And it is obvious that
$$\{\lambda:\ |i\langle k,\omega\rangle+\langle l,\Omega\rangle|=0,
\ 0<|k|\leq K,\ |l|\leq 2 \}$$
contains at most two points, that is, A2)' holds. So by the
non-singularity of $A$ and Corollary 1.1 2) we get the persistence
of invariant 2-tori on $S_{1}$ with unchanged tangent frequencies.

2) We consider the persistence of invariant tori on the parabola:
$$S_{2}:\ y_{1}(\lambda)=a_{1}\lambda,\ y_{2}(\lambda)=a_{2}\lambda^{2},
\ \lambda\in[1,2].$$
It is easy to verify that A1)' holds if and only if
$a_{1}a_{2}\neq 0$. And similar to 1), we obtain that A2)' holds
on $S_{2}$. So we get the same result on $S_{2}$ as in 1).

Since the eigenvalues of $JM$ are pure imaginary, the persistent
tori are elliptic. \\
{\large\bf Example 4.3.} We consider the following unperturbed
system:
$$N(y,u)=\frac{1}{2}y_{1}^{2}+\frac{1}{2}y_{2}^{2}+\frac{1}{2}y_{3}^{2}
+\frac{1}{2}(u_{1}^{2}-v_{1}^{2})
+\frac{1}{2}u_{2}^{2}-\frac{3}{2}v_{2}^{2}$$ with
\begin{align*}
\omega={ \left(
\begin{array}{c}
   y_{1}\\
   y_{2}\\
   y_{3}
\end{array}
\right) }, A={ \left(
\begin{array}{ccc}
   1 & 0 & 0\\
   0 & 1 & 0\\
   0 & 0 & 1
\end{array}
\right) }, R={ \left(
\begin{array}{cccc}
   y_{1} & 1 & 0 & 0\\
   y_{2} & 0 & 1 & 0\\
   y_{3} & 0 & 0 & 0
\end{array}
\right) },
\\
\Omega={ \left(
\begin{array}{c}
   1\\
   -1\\
   \sqrt[]{3}\\
    -\sqrt[]{3}
\end{array}
\right) }, M={ \left(
\begin{array}{cccc}
   1 & 0 & 0 & 0\\
   0 & -1 & 0 & 0\\
   0 & 0 & 1 & 0\\
   0 & 0 & 0 & -3
\end{array}
\right) },
\end{align*}
where $R$ is obtained on the hyperplane $S$: $y_{3}=a,\ a\neq 0$,
i.e. the sub-manifold we will consider. It is easy to verify that
A1)' holds on $S$. And we easily observe that the set
$$\{\lambda:|i\langle k,\omega\rangle+\langle l,\Omega\rangle|=0,
\ 0<|k|\leq K,\ |l|\leq 2\}$$
is a straight line in $S$, or empty set, that is, A2)' also holds
on $S$.

Since $A$ is always nonsingular on $S$, by Corollary 1.1 2) we
obtain the persistence of invariant 3-tori with the same tangent
frequencies as the unperturbed system. Besides, since all
eigenvalues of $JM$ are real, the surviving tori are hyperbolic.

\vskip 5mm \noindent{\bf\Large Acknowledgements}

The author expresses his sincere thanks to Professor Yong, Li for
his instructions and many invaluable suggestions. The author is
grateful to Dr Qingdao, Huang and Baifeng, Liu for helpful
discussions.

 \vskip 5mm \noindent{\Large\bf Appendix A}

\vskip 3mm \noindent{\large\bf Lemma A.1} {\it Suppose that $g(x)$
is a m-th differentiable function on the closure $\bar{I}$ of $I$,
where $I\subset R^{1}$ is an interval. Let
$I_{h}=\{x:|g(x)|<h,x\in I\},h>0.$ If on $I$, $|g^{(m)}(x)|\geq
d>0,$ where $d$ is a
constant, then $|I_{h}|\leq ch^{\frac{1}{m}}.$ }\\
{\bf Proof.} See Lemma 2.1 in \cite{xy}.\quad$\Box$\\
 {\large\bf
Lemma A.2} {\it Let $A,B,C$ be $n\times n,m\times m,n\times m$
matrices respectively. Then the equation
$$AX+XB=C,$$
where $X$ is an $n\times m$ unknown matrix, is solvable if and
only if $I_{m}\otimes A^{\top}+B\otimes I_{n}$ is nonsingular.
Moreover
$$X=(I_{m}\otimes A^{\top}+B\otimes I_{n})^{-1}C.$$
}\\
{\bf Proof.} See Appendix in \cite{y}.\quad$\Box$ \\
{\large\bf
Lemma A.3} {\it The eigenvalues of $i\langle
k,\omega(\lambda)\rangle I_{2m}-JM,\ i\langle
k,\omega(\lambda)\rangle I_{4m^{2}}-I_{2m}\otimes (JM)-(JM)\otimes
I_{2m}$ are $i\langle k,\omega\rangle-\Omega_{j},\ i\langle
k,\omega\rangle-\Omega_{j}-\Omega_{k}, \ j,k=1,\cdots,2m,$
respectively. }\\
 {\bf Proof.} See Appendix in
\cite{y}.\quad$\Box$

\end{document}